\documentclass{amsart}

\usepackage{amsmath}
\usepackage{amssymb}
\usepackage{amsfonts}
\usepackage{amsthm}
\usepackage{mathabx}
\usepackage{mathrsfs}
\usepackage{dirtytalk}

\newcommand*\norm[1]{ \left\| #1 \right\| }

\newcommand*\Z{ \mathbb{Z} }

\newcommand*\R{ \mathbb{R} }
\newcommand*\C{ \mathbb{C} }
\newcommand*\N{ \mathbb{N} }

\newcommand*\T{\mathbb T }
\newcommand\set[1]{\left\{ #1 \right\}}

\newtheorem{mylemma}{Lemma}

\newtheorem{mythm}{Theorem}

\newtheorem{Remark}{Remark}

\title{Large deviation estimates for quasi-periodic Gevrey cocycles}
\author{Matthew Powell}
\address{Department of Mathematics, Georgia Institute of Technology, Atlanta GA, 30332}

\date{\today}

\begin{document}
\maketitle

\begin{abstract}
In this note we use an approximation scheme to establish large deviations for quasi-periodic Gevrey cocycles. As an application, we obtain continuity in the cocycle for the Lyapunov exponent.
\end{abstract}

The purpose of this note is to present a method for obtaining large deviations, and thus continuity, for Lyapunov exponents for quasi-periodic cocycles with Gevrey entries.

We say $f: \R^n \to \C$ is $s$-Gevrey (denoted $G^s(\R^n)$ if, for all multi-indexes $\alpha \in \N^n,$ there exists $C = C(f)$ such that
$$\norm{D^\alpha f}_\infty \leq C^{|\alpha|+1} |\alpha!|^s.$$
When $s = 1,$ these are the analytic funtions. Since $G^s, s > 1$ contains non-analytic functions, we cannot use analytic methods to study them directly.

An $s$-Gevrey quasiperiodic cocycle is a function $A: \T^d \to SL(2,\C)$ whose matrix elements are in the class $G^s.$ For a fixed frequency $\omega \in \T^d,$ we define the cocycle iterates via
$$A_N(x) := \prod_{j = N - 1}^0 A(x + j \omega).$$
We may define the Lyapunov exponent of such a cocycle in the usual way:
$$L_N(A, \omega) := \lim_{N\to\infty} \frac 1 N \int_{\T^d} \ln\norm{A_N(x)} dx.$$

Our goal is to understand the continuity of $L_N(A,\omega)$ in $A$ and/or $\omega.$ Following ideas introduced by Bourgain and his coauthors, the first major step to such an understanding is a large deviation estimate. Since the argument used by Bourgain explicitly uses properties of {\it analytic} cocycles, we will need to adjust our approach slightly.

Let us recall the large deviation for analytic cocycles due to Bourgain \cite{BourgainContinuity} and extended to the case of general $2\times 2$ matrices (c.f. \cite{DuarteKleinBook} and \cite{PowellContinuity}). For simplicity, we will formulate the result for invertible matrices only, since that is the setting of this note.

\begin{mythm}[Theorem 3.1 from \cite{PowellContinuity}] 
Let $A$ be an analytic $GL(2,\C)$ cocycle such that $L_N(A,x,\omega) := \frac 1 N \ln\norm{A_N(x)}$ has a 1-bounded plurisubharmonic extension to the complex strip $\set{z \in \C^d: |\Im z_j| < \rho}.$ Moreover, suppose $\omega \in \T^d$ is such that
$$\norm{k\cdot \omega} > \delta$$
for all $0 < |k| \leq K_0$ and 
$$N \geq K_0 \delta^{-1}.$$
Then $L_N(A,x,\omega)$ satisfies a large deviation estimate: for some $c = c(d) \in (0, 1]$
$$\left|\set{x\in \T^d: \left|L_N(A,x,\omega) - L_N(A,\omega)\right| > \rho^{-1} K_0^{-c}}\right| \leq e^{-\rho K_0^c}.$$
\end{mythm}

Our goal here is to establish a suitably weak analogue of the above estimate for quasi-periodic Gevrey cocycles. In fact, we will prove the following quantitative large deviation estimate.

\begin{mythm}[Locally uniform Large Deviation for Diophantine frequency and restricted Gevrey class]\label{thm:GLDT}
Let $A$ be an $s$-Gevrey $SL(2,\C)$ cocycle with $s < 1 + \frac{c(d)}{1 + b}$ and suppose $\omega \in \T^d$ is such that
$$\norm{k\cdot \omega} > \kappa |k|^{-b}$$
for all $0 < |k| \leq K_0$. Then for 
$$N \sim K_0^{1 + b}$$
the finite-scale Lyapunov exponent, $L_N(A,x,\omega)$, satisfies a large deviation estimate: for some $c = c(d) \in (0, 1]$ and $C = C(A) <\infty,$
$$\left|\set{x\in \T^d: \left|L_N(A,x,\omega) - L_N(A,\omega)\right| > 2C^{1/s} N^{(s - 1)-\frac{c}{1 + b}}}\right| \leq e^{-\frac 12 C^{-1/s} N^{(1 - s) + \frac{c}{1 + b}}}.$$
Moreover, the constant $C(A)$ depends continuously on the cocycle $A.$
\end{mythm}

\begin{Remark}
Note that when $s = 1$ this is precisely the analytic version of the LDT.
\end{Remark}
\begin{Remark}
We can replace $SL(2,\C)$ with $GL(2,\C)$ in this theorem with a small technical argument.
\end{Remark}
\begin{Remark}
Since we are restricting our attention to $s$-Gevrey $SL(2\C)$ cocycles, we lose nothing by using $\log_+ = \max(\log, 0)$ in place of $\log.$ Throughout the rest of this paper, $\log$ will mean $\log_+.$ We make this choice to avoid a technical argument needed to control $\|A_N^{m)}\|$ from below.
\end{Remark}

In order to prove Theorem \ref{thm:GLDT}, 
we will argue via approximation. The idea of approximating Gevrey functions by suitable analytic functions is fairly standard, and has been used by many authors to extend analytic results to the Gevrey classes (c.f. \cite{SKleinMultiD} and \cite{SKleinOneD}). Our approach shares some thematic similarities with a method used by S. Klein to obtain localization \cite{SKleinMultiD}, but we will need to make substantial changes to overcome a technical obstacle which will allow us to consider a wider class of cocycles, and reduce the assumptions. 

\section{An approximation scheme}

In brief, let us denote by $A^{(m)}$ the analytic cocycle with entries given by a truncated Fourier series
$$(A^{(m)})_{ij}(x) = \sum_{|k|\leq m} \widehat{A_{ij}}(k) e^{i k\cdot x}.$$
Then for any $\rho = \rho(m)$ such that $L_N(A^{(m)}, x, \omega)$ has a 1-bounded plurisubharmonic extension to the complex strip $\C_\rho^d,$ the analytic large deviation estimate holds. By a standard telescoping argument, we have
$$\left|L_N(A, x, \omega) - L_N(A^{(m)}, x, \omega)\right| \leq \frac 1 N C^N \norm{A - A^{(m)}}_\infty,$$
where $C = C(A).$ Combining this estimate with the analytic large deviation, we obtain, under the same assumptions as the analytic large deviation estimate, and using $\rho = \rho(m),$
$$\left|\set{x\in \T^d: \left|L_N(A,x,\omega) - L_N(A,\omega)\right| > \rho^{-1} K_0^{-c} - 2 N C^N \norm{A - A^{(m)}}_\infty}\right| \leq e^{-\rho K_0^c}.$$
This would be a ``good'' estimate if 
$$\rho^{-1} K_0^{-c} - 2 N C^N \norm{A - A^{(m)}}_\infty \geq \frac 12 \rho^{-1} K_0^{-c},$$
which will require us to find optimal choices of $m$ and $\rho(m).$ These choices will depend on $N,$ which is the source of the major difficulty we will face.

\subsection{An optimal choice of $\rho$}
We will begin by finding an allowable range of $\rho$ which ensures that the 1-boundedness condition holds for $A^{(m)}.$ Let $g \in G^s(\T^d)$ such that $|g(x)| \geq 1$ (as a stand-in for the $SL(2,\C)$ condition), and set $g_m(x) = \sum_{|k| \leq m} \hat g(k) e^{ik\cdot x}.$ Since $G^s$ is characterized by a bound on the growth of the derivatives, we can use the Fourier transform to generate an equivalent characterization in terms of Fourier coefficients: $f \in G^s(\T^d)$ if, for some $C = C(f),$
$$\hat{f}(k) \leq e^{-C^{-1/s} |k|^{1/s}}.$$
Thus
$$\left| g(x) - g_m(x)\right| \leq e^{-C^{-1/s} m^{1/s}}.$$
Moreover, $g_m$ is entire, so we can estimate $u_m(x) = \frac{1}{N} \ln|\prod_{j = 1}^N g_m(z + j\omega)|$ on $\C^d$ via
\begin{align}
|u_m(z)| &\leq \log\left|\sum_{|k| \leq m} \hat g(k) e^{ik\cdot z}\right|\\
&\leq \log\left(\sum_{|k| \leq m} |\hat g(k)| e^{|k\cdot \Im z|}\right)\\
&\leq \log\left(\sum_{|k| \leq m}  e^{|k| |\Im z| -C^{-1/s} |k|^{1/s}}\right).
\end{align}
In order for this to be 1-bounded, we see that it is sufficient for $|\Im z|$ to satisfy
$$|\Im z| \leq \frac 12 C^{-1/s} m^{1/s - 1}.$$
Hence, we will take
$$\rho(m) \leq \frac 12 C^{-1/s} m^{1/s - 1}.$$

\subsection{An optimal choice of $m$}
Now we turn our attention to the degree of our approximation, $m.$ With $\rho(m)$ in hand, it remains to find $m$ such that
$$\rho^{-1} K_0^{-c} - 2 N C^N \norm{A - A^{(m)}}_\infty \geq \frac12\rho^{-1} K_0^{-c}.$$
Observe that 
$$\norm{A - A^{(m)}}_\infty \leq e^{-C^{-1/s} m^{1/s}}.$$
From this, we see that taking $m \geq C N^s,$ $N > N_0(C)$ sufficiently large, and $C = C(A)$ some sufficiently large constant is sufficient.

\subsection{An initial LDT}
Combining our choice of $\rho(m)$ and $m,$ we obtain the following large deviation estimate for $G^s$ cocycles.
\begin{mylemma}[Large deviation for generic frequency and arbitrary Gevrey class]
Let $A$ be an $s$-Gevrey $SL(2,\C)$ cocycle. Moreover, suppose $\omega \in \T^d$ is such that
$$\norm{k\cdot \omega} > \delta$$
for all $0 < |k| \leq K_0$ and 
$$N \geq K_0 \delta^{-1}.$$
Then $L_N(A,x,\omega)$ satisfies a large deviation estimate: for some $c = c(d) \in (0, 1]$
$$\left|\set{x\in \T^d: \left|L_N(A,x,\omega) - L_N(A,\omega)\right| > 2C^{1/s} N^{s - 1} K_0^{-c}}\right| \leq e^{-\frac 12 C^{-1/s} N^{1 - s} K_0^c}.$$
\end{mylemma}
Note that when $s = 1$ this is precisely the analytic version of the LDT. One immediately sees a problem, however, when $s > 1.$ In this case, $N^{s - 1} K_0^{-c}$ will typically be large, which means the estimate is quite trivial and of no practical use.

How might we rescue this estimate? S. Klein used an approach due to Bourgain, Goldstein, and Schlag to show that, under a strong Diophantine assumption, it is possible to improve this estimate inductively for Schr\"odinger cocycles with sufficiently large coupling constant and $G^s$ potentials which satisfy a transversality condition. This approach is restricted to Schr\"odinger-like cocycles due to the coupling constant requirement, and the use of a transversality condition is, in some sense, unnatural as a prerequisite for a large deviation estimate. That approach is, however, very apt if one wishes to obtain localization as in \cite{SKleinMultiD} where transversality is the Gevrey analogue of non-constant for analytic functions, which is essential for localization. In brief, the large coupling and transversality condition in \cite{SKleinMultiD} allowed Klein to establish an effective large deviation of the above type at some fixed initial scale $N_0.$ These properties then allowed Klein to run an induction scheme in the spirit of \cite{BourgainGoldsteinSchlag} to boost the large deviation to arbitrary scale. 

We propose a different approach. Instead of restricting our attention to Schr\"odinger cocycles and appealing to techniques from \cite{BourgainGoldsteinSchlag}, we will instead restrict the allowable range of the Gevrey exponent $s$, which will allow us to avoid both the induction scheme and transversality condition.


Let us observe that, if $\omega$ satisfies a Diophantine assumption of the form
$$\norm{k\cdot \omega} \geq \kappa |k|^{-b}; \quad 0 < |k|$$ then we may rewrite this estimate for $K_0 \sim N^{1/(1 + b)}:$
$$\left|\set{x\in \T^d: \left|L_N(A,x,\omega) - L_N(A,\omega)\right| > 2C^{1/s} N^{s - 1} N^{-c/(1 + b)}}\right| \leq e^{-\frac 12 C^{-1/s} N^{1 - s} N^{c/(1 + b)}}.$$
It is now clear that taking $s = 1 + \epsilon$ for any $\epsilon < c/(1 + b)$ ensures that this is a good enough large deviation for us to work with. Hence we arrive at the following more useful estimate for so-called restricted Diophantine frequencies (that is, frequencies which satisfy a Diophantine condition up to some finite scale $K_0$).

\begin{mylemma}[Locally uniform Large Deviation for (restricted) Diophantine frequency and restricted Gevrey class]\label{lem:GLDT}
Let $A$ be an $s$-Gevrey $SL(2,\C)$ cocycle with $s < 1 + \frac{c(d)}{1 + b}$ and suppose $\omega \in \T^d$ is such that
$$\norm{k\cdot \omega} > \kappa |k|^{-b}$$
for all $0 < |k| \leq K_0$. Then for 
$$N \sim K_0^{1 + b}$$
the finite-scale Lyapunov exponent, $L_N(A,x,\omega)$, satisfies a large deviation estimate: for some $c = c(d) \in (0, 1]$ and $C = C(A) <\infty,$
$$\left|\set{x\in \T^d: \left|L_N(A,x,\omega) - L_N(A,\omega)\right| > 2C^{1/s} N^{(s - 1)-\frac{c}{1 + b}}}\right| \leq e^{-\frac 12 C^{-1/s} N^{(1 - s) + \frac{c}{1 + b}}}.$$
Moreover, the constant $C(A)$ depends continuously on the cocycle $A.$
\end{mylemma}

Theorem \ref{thm:GLDT} now follows immediately from Lemma \ref{lem:GLDT} with $\omega$ satisfying an unrestricted Diophantine condition: $\norm{k\cdot \omega} > \kappa |k|^{-b}$ for all $k \in \Z^d\backslash\{0\}.$

\section{Application}
It was shown by Bourgain \cite{BourgainContinuity} (see also \cite{DuarteKleinBook} and \cite{PowellContinuity}) that, once one has a locally uniform large deviation of the above type for Diophantine frequencies, one also has continuity of the associated Lyapunov exponent. We thus obtain the following as an immediate corollary of the above.

\begin{mythm}
Suppose $\omega \in \T^d$ satisfies 
$$\norm{k\cdot \omega} \geq C |k|^{-b}.$$
Then $L(A,\omega)$ is a continuous function of $A$ in the $G^s(\T^d)$ topology for any $1 \leq s< 1 + c/(1 + b).$ 
\end{mythm}

\bibliographystyle{abbrv} 
\bibliography{GevreyCocycles}

\end{document}